\documentclass{amsart}

\usepackage{graphicx}
\usepackage{amssymb}
\newtheorem{theorem}{Theorem}[section]

\newtheorem{corollary}[theorem]{Corollary}

\theoremstyle{definition}

\newtheorem{example}[theorem]{Example}

\newcommand{\C}{\Bbb C}
\newcommand{\Z}{\Bbb Z}
\newcommand{\D}{\Delta}

\newcommand{\p}{\partial}
\newcommand{\la}{\langle}
\newcommand{\ra}{\rangle}

\theoremstyle{remark}
\newtheorem{remark}[theorem]{Remark}

\numberwithin{equation}{section}

\begin{document}

\title
{A Torres condition for twisted Alexander polynomials}

\author{Takayuki Morifuji}

\thanks{2000 {\it Mathematics Subject Classification}. 
Primary 57M25, Secondary 57M05, 57M27}

\thanks{{\it Key words and phrases.\/} 
Torres condition, link group, twisted Alexander polynomial.}

\address{Department of Mathematics, 
Tokyo University of Agriculture and Technology, 
Koganei, Tokyo 184-8588, Japan}

\email{morifuji@cc.tuat.ac.jp}

\begin{abstract}
As a generalization of a fundamental result 
about the Alexander polynomial of links, 
we give a description of a Torres condition for 
the twisted Alexander polynomial of links 
associated to a unimodular representation. 
\end{abstract}

\maketitle

\section{Introduction}

The theory of twisted Alexander polynomial 
was introduced by Lin \cite{Lin01-1} and 
Wada \cite{Wada94-1} independently. 
Lin defined it for knots in the $3$-sphere 
using regular Seifert surfaces. 
On the other hand, 
Wada defined the twisted Alexander polynomial for 
finitely presentable groups, 
which include the link groups. 
In particular, 
as an application, 
Wada told the Kinoshita-Terasaka knot 
from the Conway knot by means of his invariant. 
Shortly afterward, 
several significant results on 
the original Alexander polynomial were 
generalized to the twisted case. 
For example, 
equivalence of the twisted Alexander polynomial 
and the Reidemeister torsion, 
and its symmetry \cite{Kitano96-1}, \cite{KL99-1}, 
sliceness obstruction for knots and 
a relation to the Casson-Gordon invariant 
\cite{KL99-1}, \cite{KL99-2}, 
monicness of the twisted Alexander polynomial 
for fibered knots \cite{Cha03-1}, \cite{GKM02-1} 
and so on. 
Recently 
the twisted Alexander polynomials 
are extensively investigated. 
See for instance 
\cite{GM03-1}, 
\cite{HP04-1}, 
\cite{JW93-1}, 
\cite{KM04-1}, 
\cite{KSW04-1}, 
\cite{Morifuji01-1} and 
\cite{Suzuki04-1}.

However, 
almost all results mentioned above are 
basically about knots in the $3$-sphere 
and 
it seems that 
there are few generalized results on links. 
The purpose of the present paper 
is to give a generalization of 
the following well-known formula for 
the Alexander polynomial of links. 

\begin{theorem}[Torres \cite{Torres53-1}]
For 
the Alexander polynomial 
$\Delta_L(t_1,\ldots,t_\mu)$ 
of a $\mu$-component link $L=L_1\cup \cdots\cup L_\mu$, 
it holds that 
$$
\displaystyle{
\Delta_L(t_1,\ldots,t_{\mu-1},1)
=
\begin{cases}
\displaystyle{\frac{t_1^{l_1}-1}{t_1-1}}\Delta_{L'}(t_1)
& if~\mu =2 \\
(t_1^{l_1}\cdots t_{\mu-1}^{l_{\mu-1}}-1)
\Delta_{L'}(t_1,\ldots,t_{\mu-1})
& if~\mu>2, 
\end{cases}}
$$
where 
$L'=L_1\cup\cdots\cup L_{\mu-1}$ 
is the link obtained 
from $L$ by removing $L_\mu$ and 
$l_i$ denotes the linking number of the components 
$L_i$ and $L_\mu$. 
\end{theorem}

More precisely, 
we give a description of a Torres condition for 
the twisted Alexander polynomial of links 
associated to a unimodular representation. 
In the next section, 
we briefly recall the definition of 
the twisted Alexander polynomial for the link group. 
The precise statement and the proof of 
the main theorem of this paper 
are given in Section 3. 

\section{Twisted Alexander polynomial for links}

Let $L=L_1\cup\cdots\cup L_\mu$ be 
a $\mu$-component link in the $3$-sphere. 
We denote the fundamental group of its exterior $E$ 
by $G(L)$. 
Namely, 
we put $G(L)=\pi_1(E)$  
and call it the link group. 
We choose and fix a Wirtinger presentation 
of $G(L)$:
$$
G(L)=
\langle x_1,\ldots,x_u~|~r_1,\ldots,r_{u-1}\rangle.
$$
Then 
the abelianization homomorphism 
$$
\alpha:G(L)\to H_1(E;\Bbb Z)\cong\Bbb Z^{\oplus\mu}
=\la t_1\ra\oplus\cdots\oplus\la t_\mu\ra
$$ 
is given by 
assigning to each generator $x_i$ 
the meridian element $t_k\in H_1(E;\Z)$ of 
the corresponding component $L_k$ of $L$. 
In this paper, 
we consider a linear representation 
$\rho:G(L)\to SL(n;F)$, 
where $F$ denotes a field. 

These maps naturally induce two ring homomorphisms 
$\tilde{\rho}: \Bbb Z[G(L)] \rightarrow M(n;F)$ 
and 
$\tilde{\alpha}:\Bbb Z[G(L)]\rightarrow 
\Z[t_1^{\pm1},\ldots,t_\mu^{\pm1}]$, 
where $\Bbb Z[G(L)]$ is the group ring of $G(L)$ 
over $\Z$ 
and 
$M(n;F)$ is the matrix algebra of degree $n$ over $F$. 
Taking the tensor of 
$\tilde{\rho}$ and $\tilde{\alpha}$, 
we obtain a ring homomorphism 
$$
\tilde{\rho}\otimes\tilde{\alpha}:
\Bbb Z[G(L)]\to M\left(n;F[t_1^{\pm1},\ldots,t_\mu^{\pm1}]\right).
$$
Let 
$F_u$ denote the free group on 
generators $x_1,\ldots,x_u$ and 
$$
\Phi:\Bbb Z[F_u]\to M\left(n;F[t_1^{\pm1},\ldots,t_\mu^{\pm1}]\right)
$$
the composite of the surjection 
$\Bbb Z[F_u]\to\Bbb Z[G(L)]$ 
induced by the presentation and the map 
$\tilde{\rho}\otimes\tilde{\alpha}$. 

Let us consider the $(u-1)\times u$ matrix $M=M(t_1,\ldots,t_\mu)$ 
whose $(i,j)$th component is the $n\times n$ matrix 
$$
\Phi\left(\frac{\partial r_i}{\partial x_j}\right)
\in M\left(n;F[t_1^{\pm1},\ldots,t_\mu^{\pm1}]\right),
$$
where 
${\partial}/{\partial x}$ 
denotes the free differential calculus. 
This matrix $M$ is called 
the Alexander matrix of $G(L)$ 
associated to the representation $\rho$. 

For 
$1\leq j\leq u$, 
let us denote by $M_j=M_j(t_1,\ldots,t_\mu)$ 
the $(u-1)\times(u-1)$ matrix obtained from $M$ 
by removing the column 
corresponding to a generator $x_j$. 
We also regard $M_j$ as 
an $n(u-1)\times n(u-1)$ matrix with coefficients in 
$F[t_1^{\pm1},\ldots,t_\mu^{\pm1}]$. 

Then 
Wada's twisted Alexander polynomial 
of a link $L$ for a representation $\rho:G(L) \to SL(n;F)$ 
is defined to be a rational function 
$$
\D_{L,\rho}(t_1,\ldots,t_\mu)
=\frac{|M_j|}{|\Phi(x_j-1)|},
$$
where 
$|M_j|$ denotes the determinant of the matrix $M_j$, 
and 
it is well-defined up to a factor 
$\pm t_1^{nk_1}\cdots t_\mu^{nk_\mu}~(k_i\in\Bbb Z)$ 
if $n$ is odd and up to only 
$t_1^{nk_1}\cdots t_\mu^{nk_\mu}$ 
if $n$ is even 
(see \cite{Wada94-1} Section 5 for details). 

\begin{remark}
In general, 
the twisted Alexander polynomial 
for finitely presentable groups is 
a rational function, 
but 
it is actually a polynomial 
for the link groups 
(see \cite{Wada94-1} Proposition 9 and 
\cite{KM04-1} Theorem 3.1). 
\end{remark} 

\section{A Torres condition}

In this section, 
we state and prove 
a generalized Torres condition 
for the twisted Alexander polynomial of links. 
An advantage of our description here 
is that we need not separate the case for $\mu=2$ 
from $\mu>2$. 
We 
first prove the theorem in the case of 
an $SL(2;F)$-representation. 
After 
reading the proof for it, 
one can easily show the similar result 
for general cases. 

\begin{theorem}
Let $L=L_1\cup \cdots\cup L_\mu$ be a $\mu$-component link 
and $L'=L_1\cup \cdots\cup L_{\mu-1}$. 
For a given representation 
$\rho':G(L')\to SL(2;F)$, it holds that 
$$
\Delta_{L,\rho}(t_1,\ldots,t_{\mu-1},1)
=\{(t_1^{l_1}\cdots t_{\mu-1}^{l_{\mu-1}})^2
-\varepsilon_{\rho'}t_1^{l_1}\cdots t_{\mu-1}^{l_{\mu-1}}+1\}
\Delta_{L',\rho'}(t_1,\ldots,t_{\mu-1}),
$$
where $\rho:G(L)\to SL(2;F)$ is the composite of 
the natural surjection $G(L)\to G(L')$ and $\rho'$, 
$l_i$ denotes the linking number 
of $L_i$ and $L_\mu$, 
and $\varepsilon_{\rho'}$is 
an element of $F$. 
\end{theorem}

\begin{proof}
For the link group $G(L)$, 
we choose a Wirtinger presentation:
$$
G(L)=
\langle x_{ij}~|~r_{kl} \rangle,
$$
where 
$x_{i1},x_{i2},\ldots,x_{ij_i}~(1\leq i\leq \mu)$ 
are generators corresponding to the component $L_i$ 
and the relation 
$$
r_{kl}=x_{k'l'}^{\pm1}x_{kl}x_{k'l'}^{\mp1}x_{k,l+1}^{-1}
$$
corresponds to a crossing of $L_{k'}$ over $L_k$.  
We should note that 
the link group $G(L)$ has the deficiency one. 

Let us consider the Alexander matrix of $G(L)$ 
associated to the representation 
$\rho:G(L)\to SL(2;F)$:
$$
M(t_1,\ldots,t_\mu)
=\left(\Phi\left(\frac{\p r_{kl}}{\p x_{ij}}\right)\right).
$$
Then 
we know that if 
we remove the column corresponding to 
a generator $x_{ij}$, 
$$
| M_{ij}(t_1,\ldots,t_\mu)|
=| \Phi(x_{ij}-1)|\D_{L,\rho}(t_1,\ldots,t_\mu)
$$
holds.
Thus 
if we set $t_\mu=1$ in $M(t_1,\ldots,t_\mu)$, 
it follows that 
$$
| M_{ij}(t_1,\ldots,t_{\mu-1},1)|
=| \Phi(x_{ij}-1)|\D_{L,\rho}(t_1,\ldots,t_{\mu-1},1)
$$
if $i\not=\mu$. 

Now 
the generators $\{x_{\mu j}\}$ appear 
in the following two kinds of relations:
$$
\mathrm{(i)}~
r_{\mu j}=x_{rs}^{\pm1}x_{\mu j}x_{rs}^{\mp1}x_{\mu,j+1}^{-1}
\quad\mathrm{and}\quad
\mathrm{(ii)}~r_{pq}=x_{\mu l}^{\pm1}x_{pq}x_{\mu l}^{\mp1}x_{p,q+1}^{-1},
$$
where 
the relation (i) corresponds to crossings of $L_r$ over $L_\mu$ 
and (ii) corresponds to 
that of $L_\mu$ over $L_p$. 
Let us see which are the contributions of these relations 
to the matrix $M(t_1,\ldots,t_{\mu-1},1)$. 
First 
the contributions of $r_{\mu j}$ are as follows:
\begin{align*}
\Phi\left(\frac{\p r_{\mu j}}{\p x_{rs}}\right)_{t_\mu=1}
&=O,\\
\Phi\left(\frac{\p r_{\mu j}}{\p x_{\mu j}}\right)_{t_\mu=1}
&=
\begin{cases}
t_r^{\pm1}\rho(x_{rs})^{\pm1} & \mathrm{if}~\mu\not= r\\
I & \mathrm{if}~\mu=r,
\end{cases}\\
\Phi\left(\frac{\p r_{\mu j}}{\p x_{\mu,j+1}}\right)_{t_\mu=1}
&=-I,
\end{align*}
where 
$O$ and $I$ denote the zero and the identity matrix respectively. 
We have used here 
the fact that 
$\rho(x_{\mu j})=I$ for $1\leq j\leq j_\mu$ 
(because the generators $\{x_{\mu j}\}$ 
are in the kernel of the natural surjective homomorphism 
from $G(L)$ to $G(L')$). 
Next 
the contributions of $r_{pq}$ are as follows:
\begin{align*}
\Phi\left(\frac{\p r_{pq}}{\p x_{\mu l}}\right)_{t_\mu=1}
&=\pm(I-t_p\rho(x_{pq})),\\
\Phi\left(\frac{\p r_{pq}}{\p x_{pq}}\right)_{t_\mu=1}
&=I,\\
\Phi\left(\frac{\p r_{pq}}{\p x_{p,q+1}}\right)_{t_\mu=1}
&=-\rho(x_{pq})\rho(x_{p,q+1})^{-1}~\mathrm{if}~p\not=\mu,
\end{align*}
and the case $p=\mu$ has already been considered. 
Therefore 
we see that the matrix $M(t_1,\ldots,t_{\mu-1},1)$ 
has the following form:
$$
M(t_1,\ldots,t_{\mu-1},1)
=
\begin{pmatrix}
A&B\\O&C
\end{pmatrix},
$$
where 
$$
A=\left(
\Phi\left(\frac{\p r_{kl}}{\p x_{ij}}\right)_{t_\mu=1}
\right)~
(k,i\not=\mu),~
B=\left(
\Phi\left(\frac{\p r_{kl}}{\p x_{\mu j}}\right)_{t_\mu=1}
\right)~
(k\not=\mu,~1\leq j \leq j_\mu),
$$
and
\begin{align*}
C
&=\left(
\Phi\left(\frac{\p r_{\mu l}}{\p x_{\mu j}}\right)_{t_\mu=1}
\right)~
(1\leq j,l\leq j_\mu)\\
&=
\begin{pmatrix}
t_{r_1}^{\delta_{r_1}}\rho(x_{r_1s_1})^{\delta_{r_1}} & -I & & & \\
&t_{r_2}^{\delta_{r_2}}\rho(x_{r_2s_2})^{\delta_{r_2}} & -I& & \\
& &t_{r_3}^{\delta_{r_3}}\rho(x_{r_3s_3})^{\delta_{r_3}} & & \\
& & &\ddots &-I \\
-I & & & & t_{r_{j_\mu}}^{\delta_{r_{j_\mu}}}
\rho(x_{r_{j_\mu}s_{j_\mu}})^{\delta_{r_{j_\mu}}}
\end{pmatrix}.
\end{align*}

In the submatrix $C$, 
there is an appearance of 
$t_i^{\delta_i}$ for each crossing of 
$L_i$ over $L_\mu~(1\leq i\leq \mu)$ 
and 
$\delta_i=1$ or $-1$ according as 
$L_i$ crosses over $L_\mu$ 
from left to right or 
from right to left. 
Thereby, 
we obtain
\begin{align*}
|C|
&= \prod_{i=1}^{j_\mu}
|\rho(x_{r_is_i})^{\delta_{r_i}}|t_{r_i}^{2\delta_{r_i}}
-\sum_{\sigma\in S}c_1^\sigma \cdots c_{j_{\mu}}^\sigma 
t_{r_1}^{\delta_{r_1}} \cdots t_{r_{j_\mu}}^{\delta_{r_{j_\mu}}}+1\\
&=(t_1^{l_1}\cdots t_{\mu-1}^{l_{\mu-1}})^2
-\varepsilon_{\rho'}t_1^{l_1}\cdots t_{\mu-1}^{l_{\mu-1}}+1,
\end{align*}
where 
$c_i^\sigma\in F$ 
denotes a component of 
the matrix $\rho(x_{r_is_i})^{\delta_{r_i}}$ 
determined by a permutation $\sigma$ 
and $S$ is a subset of the symmetric group $\mathfrak{S}_{2j_\mu}$ 
consisting of permutations which 
choose just one component from each submatrix 
$\rho(x_{r_is_i})^{\delta_{r_i}}$. 
Furthermore 
the submatrix $A$ is equivalent to 
the Alexander matrix 
$M'(t_1,\ldots,t_{\mu-1})$ of $G(L')$ 
associated to the representation 
$\rho':G(L')\to SL(2;F)$. 
Hence 
if we remove a column corresponding to 
a generator $x_{ij}~(i\not=\mu)$, 
then 
we have 
\begin{align*}
|M_{ij}(t_1,\ldots,t_{\mu-1},1)|
&= | A_{ij}||C| \\
&=\{(t_1^{l_1}\cdots t_{\mu-1}^{l_{\mu-1}})^2
-\varepsilon_{\rho'}t_1^{l_1}\cdots t_{\mu-1}^{l_{\mu-1}}+1\}
|M_{ij}'(t_1,\ldots,t_{\mu-1})|,
\end{align*}
where 
$A_{ij}$ is the matrix obtained from 
$A$ by removing the column 
corresponding to $x_{ij}$. 
Therefore, 
by definition of the twisted Alexander polynomial, 
we see that 
$$
\Delta_{L,\rho}(t_1,\ldots,t_{\mu-1},1)
=\{(t_1^{l_1}\cdots t_{\mu-1}^{l_{\mu-1}})^2
-\varepsilon_{\rho'}t_1^{l_1}\cdots t_{\mu-1}^{l_{\mu-1}}+1\}
\Delta_{L',\rho'}(t_1,\ldots,t_{\mu-1}).
$$
This completes the proof of Theorem 3.1.
\end{proof}

\begin{remark}
The fact that 
$\D_{L,\rho}(t_1,\ldots,t_{\mu-1},1)$ 
is divisible by $\D_{L',\rho'}(t_1,\ldots,t_{\mu-1})$ 
also follows from a recent result of 
Kitano, Suzuki and Wada in  \cite{KSW04-1}. 
However, 
we can have no detailed information 
on the quotient from their result. 
\end{remark}

A linear representation 
$\rho:G(L)\to GL(n;F)$ is called reducible 
if it has a nontrivial invariant subspace in $F^n$. 
In this case, 
we can obtain a piece of information 
about the coefficient $\varepsilon_{\rho'}$. 

\begin{corollary}
Under the setting as in Theorem 3.1, 
if $\rho':G(L')\to SL(2;F)$ is a reducible representation, 
then 
we have 
$$
\varepsilon_{\rho'}
=\lambda^l+\lambda^{-l}\quad
(l=l_1+\cdots+l_{\mu-1}),
$$
where 
$\lambda$ is an eigenvalue of the image of 
a generator of $G(L')$.
\end{corollary}

\begin{proof}
First 
we can assume that 
the images of generators in a Wirtinger presentation 
of $G(L)$ have the following forms:
$$
\rho(x_{ij})
=\begin{pmatrix}
a_{ij}& b_{ij}\\
0 & a_{ij}^{-1}
\end{pmatrix}~
(i\not=\mu)\quad
\mathrm{and}\quad
\rho(x_{\mu j})=I,
$$
where 
$a_{ij}\in F^\times$ and $b_{ij}\in F$. 
Because 
the representation $\rho'$ has 
a $1$-dimensional invariant subspace in $F^2$. 

Since 
$x_{ij}x_{kl}^{-1}~(i,k\not=\mu)$ is an element of 
the commutator subgroup $[G(L),G(L)]$, 
we see that $a_{ij}=a_{kl}$ holds for 
these generators. 
We then put 
$\lambda=a_{ij}$ for simplicity. 
Each lower left component of 
$\rho(x_{ij})$ is zero, 
so that 
the nontrivial terms appeared in the coefficient of 
$t_{r_1}^{\delta_{r_1}}\cdots t_{r_{j_\mu}}^{\delta_{r_{j_\mu}}}$ 
are just 
$$
-\lambda^{\delta_{r_1}+\cdots+\delta_{r_{j_\mu}}}
-\lambda^{-(\delta_{r_1}+\cdots+\delta_{r_{j_\mu}})}
=-(\lambda^l+\lambda^{-l}),
$$
where 
$l=\sum l_i$. 
This completes the proof.
\end{proof}

\begin{example}
Let 
$\rho':G(L')\to SL(2;F)$ be 
a reducible representation 
of the knot $L'=L_1$. 
Then 
the twisted Alexander polynomial of $L'$ 
associated to $\rho'$ is given by 
$$
\D_{L',\rho'}(t_1)
=\frac{\D_{L'}(\lambda t_1)\D_{L'}(\lambda^{-1}t_1)}
{(t_1-\lambda)(t_1-\lambda^{-1})},
$$
where 
$\D_{L'}(t_1)$ is the original Alexander polynomial of $L'$ 
(see the proof of \cite{KM04-1} Theorem 3.1 for instance). 
Hence 
we have 
\begin{align*}
\D_{L,\rho}(1,1)
&=\frac{(1-\lambda^{l_1})(1-\lambda^{-l_1})}{(1-\lambda)(1-\lambda^{-1})}
\D_{L'}(\lambda)\D_{L'}(\lambda^{-1})\\
&=(1+\lambda+\cdots+\lambda^{l_1-1})
(1+\lambda^{-1}+\cdots+\lambda^{-(l_1-1)})
\D_{L'}(\lambda)\D_{L'}(\lambda^{-1}).
\end{align*}
In particular, 
if $\rho'$ is trivial 
(namely, $\lambda=1$), 
then 
we obtain 
$\D_{L,\rho}(1,1)={l_1}^2$ 
(because $\D_{L'}(1)=\pm1$). 
\end{example}

\begin{example}
Let 
$\rho':G(L')\to SL(2;F)$ be the trivial representation. 
In this case 
$\varepsilon_{\rho'}=2$ holds, 
so that 
we have
$$
\Delta_{L,\rho}(t_1,\ldots,t_{\mu-1},1)
=(t_1^{l_1}\cdots t_{\mu-1}^{l_{\mu-1}}-1)^2
\Delta_{L',\rho'}(t_1,\ldots,t_{\mu-1}).
$$
This formula corresponds to the square of Torres' 
original formula in Theorem 1.1. 
In particular, 
$\D_{L,\rho}(1,\ldots,1)=0$ holds for $\mu> 2$.
\end{example}

If 
we slightly modify 
the proof of Theorem 3.1, 
we obtain the following general formula for 
a unimodular representation 
$\rho':G(L')\to SL(n;F)$. 
We omit here 
the repetitious proof. 

\begin{theorem}
Let $L=L_1\cup \cdots\cup L_\mu$ be a $\mu$-component link 
and $L'=L_1\cup \cdots\cup L_{\mu-1}$. 
For a given representation 
$\rho':G(L')\to SL(n;F)$, it holds that 
\begin{align*}
\Delta_{L,\rho}(t_1,\ldots,t_{\mu-1},1)
&=\{(t_1^{l_1}\cdots t_{\mu-1}^{l_{\mu-1}})^n
+\sum_{k=1}^{n-1}
\varepsilon_{k,\rho'}(t_1^{l_1}\cdots t_{\mu-1}^{l_{\mu-1}})^{n-k}
+(-1)^n\}\\
&\quad \times \Delta_{L',\rho'}(t_1,\ldots,t_{\mu-1}),
\end{align*}
where $\rho:G(L)\to SL(n;F)$ is the composite of 
the natural surjection $G(L)\to G(L')$ and $\rho'$, 
$l_i$ denotes the linking number 
of $L_i$ and $L_\mu$, 
and 
$\varepsilon_{k,\rho'}~(1\leq k\leq n-1)$ are 
elements of $F$. 
\end{theorem}

Finally, 
we extend Corollary 3.3 
when 
all the images of the representation 
$\rho':G(L')\to SL(n;F)$ are 
upper triangle matrices. 

\begin{corollary}
Under 
the setting as in Theorem 3.6, 
if\/ $\mathrm{Im}(\rho')$ are upper triangle matrices, 
then 
the coefficient $\varepsilon_{k,\rho'}$ 
is given by 
$$
\varepsilon_{k,\rho'}
=(-1)^k\sum_{1\leq i_1<\cdots<i_k\leq n}
(\lambda_1\cdots\hat{\lambda}_{i_1}\cdots\hat{\lambda}_{i_k}
\cdots\lambda_n)^l,
$$
where 
$\lambda_k~(1\leq k\leq n)$ 
are the eigenvalues of the image of 
a generator of $G(L')$ and 
$\hat{\lambda}_k$ implies that $\lambda_k$ 
is removed from the product. 
\end{corollary}

\noindent
\textit{Acknowledgements}. 
This paper was written 
while the author was visiting 
the Ludwig-Maximilians-Universit\"{a}t in M\"{u}nchen. 
He would like to express his sincere thanks for their hospitality.

\bibliographystyle{amsplain}

\begin{thebibliography}{W}

\bibitem{Cha03-1}
J. C. Cha, 
\textit{Fibred knots and twisted Alexander invariants}, 
Trans. Amer. Math. Soc. {\bf 355} (2003), 4187--4200.


\bibitem{GKM02-1}
H. Goda, T. Kitano and T. Morifuji, 
\textit{Reidemeister torsion, twisted Alexander polynomial 
and fibered knots}, Comment. Math. Helv. {\bf 80} (2005), 
51--61.

\bibitem{GM03-1}
H. Goda and T. Morifuji, 
\textit{Twisted Alexander polynomial for $SL(2,\C)$-representations 
and fibered knots}, 
C. R. Math. Acad. Sci. Soc. R. Can. {\bf 25} (2003), 97--101.

\bibitem{HP04-1}
M. Heusener and J. Porti, 
\textit{Deformations of reducible representations 
of $3$-manifold groups into $PSL(2,\C)$}, 
math.GT/0411365. 

\bibitem{JW93-1}
B. Jiang and S. Wang, 
\textit{Twisted topological invariants associated with representations}, 
in Topics in Knot Theory (1993), 211--227.

\bibitem{KL99-1}
P. Kirk and C. Livingston, 
\textit{Twisted Alexander invariants, Reidemeister torsion, and Casson-Gordon 
invariants}, Topology {\bf 38} (1999), 635--661.

\bibitem{KL99-2}
P. Kirk and C. Livingston, 
\textit{Twisted knot polynomials: inversion, mutation and concordance}, 
Topology {\bf 38} (1999), 663--671.

\bibitem{Kitano96-1}
T. Kitano, 
\textit{Twisted Alexander polynomial and Reidemeister torsion}, 
Pacific J. Math. {\bf 174} (1996), 431--442.

\bibitem{KM04-1}
T. Kitano and T. Morifuji, 
\textit{Divisibility of twisted Alexander polynomial 
and fibered knots}, preprint. 

\bibitem{KSW04-1}
T. Kitano, M. Suzuki and M. Wada, 
\textit{Twisted Alexander polynomial and surjectivity 
of a group homomorphism\/}, preprint.

\bibitem{Lin01-1}
X. S. Lin, 
\textit{Representations of knot groups and twisted Alexander polynomials}, 
Acta Math. Sin.
(Engl. Ser.) {\bf 17} (2001), 361--380.

\bibitem{Morifuji01-1}
T. Morifuji, 
\textit{Twisted Alexander polynomial for the braid group}, 
Bull. Austral. Math. Soc. {\bf 64} (2001), 1--13.

\bibitem{Suzuki04-1}
M. Suzuki, 
\textit{Twisted Alexander polynomial for the Lawrence-Krammer 
representation\/}, 
Bull. Austral. Math. Soc. {\bf 70} (2004), 67--71. 

\bibitem{Torres53-1} 
G. Torres, 
\textit{On the Alexander polynomial}, 
Ann. of Math. {\bf 57} (1953), 57--89.

\bibitem{Wada94-1} M. Wada, 
\textit{Twisted Alexander polynomial for finitely 
presentable groups},
Topology {\bf 33} (1994), 241--256.

\end{thebibliography}

\end{document}